\documentclass[Journal,9pt]{IEEEtran}
\usepackage{amssymb}
\usepackage{graphicx,subfigure}
\usepackage{latexsym,amsmath,amsthm,amssymb,amsfonts}
\usepackage{mathrsfs}
\usepackage{amscd}
\usepackage[ansinew]{inputenc}
\usepackage{color}
\usepackage{cite}
\usepackage[all]{xy}
\usepackage{balance}
\newtheorem{theorem}{Theorem}
\newtheorem{lemma}{Lemma}

\def\QED{\mbox{\rule[0pt]{1.5ex}{1.5ex}}}
\def\endproof{\hspace*{\fill}~\QED\par\endtrivlist\unskip}

\newtheorem{assumption}{Assumption}

\newtheorem{remark}{Remark}

\newcommand{\OMIT}[1]{}

\renewcommand{\baselinestretch}{1.1}

\begin{document}


\title{Distributed optimization with nonconvex velocity constraints, nonuniform position constraints and nonuniform stepsizes\thanks{Peng Lin, Chunhua Yang and Weihua Gui are with the School of Information Science and Engineering, Central South University, Changsha, China. Wei Ren is with the Department of Electrical and Computer Engineering, University of California, Riverside, USA.
E-mail: lin$\_$peng0103@sohu.com,
ren@ee.ucr.edu, ychh@csu.edu.cn, gwh@csu.edu.cn.}} 
\author{\mbox{Peng Lin}, Wei Ren, Chunhua Yang and Weihua Gui}                                           


       \markboth{IEEE Transactions on Automatic Control.}
        {IEEE Transactions on Automatic Control.}

\maketitle

\begin{abstract}                          
This note is devoted to the distributed optimization problem of multi-agent systems with nonconvex velocity constraints, nonuniform position constraints and
nonuniform stepsizes. Two distributed constrained algorithms with nonconvex velocity constraints and nonuniform stepsizes are proposed in the absence and the presence of nonuniform position constraints by introducing a switching mechanism to guarantee all agents' position states to remain in a bounded region. The algorithm gains need not to be predesigned and can be selected by each agent using its own and neighbours' information. By a model transformation, the original nonlinear time-varying system is converted into a linear time-varying one with a nonlinear error term. Based on the properties of stochastic matrices, it is shown that the optimization problem can be solved as long as the communication topologies are jointly strongly connected and balanced.
Numerical examples are given to show the obtained theoretical results.

\end{abstract}


\section{Introduction}
As an important branch of distributed control theory, the distributed optimization problem of multi-agent systems has attracted
more and more attention from the control community \cite{angelia, Elia, cotes, Cortes3, shi, liu, yuan, yuan1, Hong, Zhu, lup, Kvaternik, Zavlanos, srivast, plinwrenysong, linren4, linren5}.
For example, {{articles}} \cite{angelia,shi,liu} studied the distributed optimization problems with and without convex constraints by a projection algorithm
 and showed that all agents reach a consensus while optimizing the given team performance functions when the communication topologies are jointly strongly connected and balanced. By introducing an integrator term in the algorithm for each agent, articles \cite{Elia,cotes,Cortes3} solved the distributed optimization problem without using the vanishing stepsizes when the communication topologies are strongly connected balanced directed graphs.

 When the constraints are taken into account, most of the existing works assumed the constraint sets to be convex and few works have paid attention to the case of nonconvex constraints. In practical applications, the constraint sets might not be convex, e.g., the velocities of the quadrotors. It is meaningful to study the distributed optimization problem with nonconvex constraint sets. Moreover, most of the existing works assumed the stepsizes of the gradients or subgradients to be uniform at any instant, which made the algorithms there not fully distributed. {Article \cite{srivast} took nonuniform stepsizes into account for the distributed optimization problem in a stochastic setting, but through computing its mathematical expectation, the proposed algorithms are essentially of the uniform-stepsize ones.} Articles \cite{linren4,linren5} introduced a kind of state-dependent stepsizes to enable each agent to be able to use its own and neighors' information to optimize the team performance function without using predesigned stepsizes, but the constraint sets were assumed to be convex and nonconvex constraint sets were not considered.

%
%


In this paper, we focus on distributed optimization of multi-agent systems with nonconvex velocity constraints, nonuniform position constraints and
nonuniform stepsizes. In \cite{linren3}, the consensus problem of multi-agent system with nonconvex velocity constraints was studied, but due to the nonlinearity of the optimization term and the unbalance of the agent interaction there caused by the nonlinear constraint operators, the results cannot be applied to solve the optimization problem. Besides, most of the existing works only considered the position constraints and none has taken into account the position and velocity constraints simultaneously. Due to the nonlinearity caused by the nonconvex velocity constraints, the nature of the system is totally changed, which makes the existing approaches no longer valid for our setting.
To solve the nonconvex velocity-constrained optimization
problem, we first propose a distributed constrained algorithm by introducing a switching mechanism to guarantee all agents' position states to remain in a bounded region. Thereinto, the algorithm gains need not to be predesigned and can be selected by each agent using its own and neighbours' information.
Second, by using a model transformation, we convert the original nonlinear time-varying system into a linear time-varying one with a nonlinear error term. Third, based on the properties of stochastic matrices, it is shown that the effects of the error term on the consensus convergence of the system vanish to zero as time evolves and the optimization problem can be solved if the communication topologies are jointly strongly connected and balanced. After that, we extend the results to the case when there further exist nonuniform convex position constraints. A distributed constrained algorithm using the projection operator in combination with the switching mechanism is proposed and it is shown that the optimization problem can be solved by combing the above analysis approach with that in \cite{plinwrenysong}. {Compared with the existing works, the main contribution of this note is that the three challenges, namely, nonconvex velocity constraints, nonuniform position constraints and nonuniform stepsizes, are addressed simultaneously for distributed optimization of multi-agent systems. In fact, each of the above challenges is rarely addressed in the literature, not to mention a combination of them.  Also the agent dynamics under consideration are in the form of double integrators instead of single integrators. In addition, the proposed algorithms are fully distributed and can be implemented by using only local information and local interaction.}


\section{Preliminaries}
In this section, some preliminary results about graph theory, projection operator and stochastic matrices are introduced (see \cite{s10}, \cite{Facchinei} and \cite{s11}).
Let $\mathcal{G}(\mathcal{I},\mathcal{E})$ be a directed
graph, where $\mathcal{I}=\{1,\cdots,n\}$ is the node set, and $\mathcal{E}\subseteq\mathcal{I}\times \mathcal{I}$ is the edge
set.
An edge of $\mathcal{G}$, denoted by $(j,i)$, represents the information flow from node $j$ to node $i$. It is assumed that $(i,i)\notin \mathcal{E}$ for all $i$.
The neighbor set of node $i$ is denoted by $\mathcal{N}_i=\{j\in
\mathcal{I}:(j,i)\in \mathcal{E}\}$. The edge weight of $(j,i)$ is defined such that  $a_{ij}>\mu_c$ for some constant $\mu_c>0$ if $(j,i) \in \mathcal{E}$ and $a_{ij}=0$ otherwise. The Laplacian of the directed graph $\mathcal{G}$, denoted by $L$,
is defined as $\lfloor {L}\rfloor_{ii} =\sum_{j=1}^na_{ij}$ and $\lfloor {L}\rfloor_{ij}=-a_{ij}$ for all $i\neq j$, where $\lfloor {L}\rfloor_{ii}$ and $\lfloor {L}\rfloor_{ij}$
denote the $ii$th and $ij$th entries of the matrix ${L}$. For a given group of nodes,
the union of a set of graphs is a graph whose edge sets are the union of the edge sets of the graphs in
the set. A directed path is a sequence of ordered edges of the
form $({i_1},{i_2}),({i_2},{i_3}), \cdots,$ where ${i_j}\in
\mathcal{I}$. A directed graph is
strongly connected if there is a directed path from every node to
every other node.

\begin{lemma}\label{le1u}{\rm \cite{Facchinei} Suppose that {$Y\neq\emptyset$} is a closed convex set in {$\mathbb{R}^m$}. The following statements hold.\\
(1) For any $y\in \mathbb{R}^m$, $\|y-P_Y(y)\|$ is continuous with respect to $y$ and $\nabla \frac{1}{2}\|y-P_Y(y)\|^2=y-P_Y(y)$ where $P_{Y}$ denotes the projection operator defined as $P_{Y}(y)=\mathrm{arg} \min\limits_{{s}\in
Y}\|y-{s}\|$;\\
(2) For any {$y,z\in\mathbb{R}^m$} and
all {$\xi\in Y$}, $[y-P_Y(y)]^T(y-\xi)\geq0$,
{$\|P_Y(y)-\xi\|^2\leq\|y-\xi\|^2-\|P_Y(y)-y\|^2$}
and {$\|P_Y(y)-P_Y(z)\|\leq\|y-z\|.$}}\end{lemma}

Given $C\in \mathbb{R}^{n\times r}$,
$C$ is nonnegative ($C\geq0$) if all its elements are
nonnegative, and  $C$ is positive ($C>0$) if all its
elements $c_{ij}$ are positive. If a nonnegative matrix $C\in
\mathbb{R}^{n\times n}$ satisfies $C\textbf{1}=\textbf{1}$, then it
is {stochastic}.

\section{Model}
Consider a multi-agent system with $n$ agents. Let $\mathcal{G}(kT)$ denote its communication graph, where $k$ is the discrete time index and $T$ is the
sampling period, ${L}(kT)$ denote the Laplacian of $\mathcal{G}(kT)$ and $\mathcal{N}_{i}(kT)$ denote the neighbor set whose information agent $i$ has access to. Suppose that each agent has the following dynamics
\begin{eqnarray}\label{eq11}\begin{array}{lll}
{r}_i((k+1)T)&=&r_i(kT)+v_i(kT)T\\
{v}_i((k+1)T)&=&u_i(kT)
\end{array}\end{eqnarray}
where $r_i(kT)\in \mathbb{R}^m$, $v_i(kT)\in \mathbb{R}^m$ and $u_i(kT)\in \mathbb{R}^m$ are
the position, velocity and
the control input of agent $i$ for some positive integer $m$. In the following, all
``$(kT)$" will be replaced by ``$(k)$" when no confusion arises. In reality, the agent velocities are often constrained to remain in nonconvex sets. For example, quadrotors can move towards every direction but the maximum velocities in different directions might be different and all of their possible velocities do not necessarily form a convex set. To this end,
we assume that
$v_i(k)\in V_i\subseteq \mathbb{R}^m$ where each $V_i$ is a nonconvex set that is known to only agent $i$. Before giving the specific assumption, we need first introduce a constraint operator that will also be used in our algorithms.

Define \begin{eqnarray*}\mathrm{S}_{V_i}(x)=\left\{\begin{array}{lll}
\frac{\displaystyle x}{\displaystyle\|x\|}
\max\limits_{{0\leq\beta\leq\|x\|}}{\Big\{}\beta{\big|}\frac{\displaystyle \alpha \beta x}{\displaystyle\|x\|}\in V_i,  \forall 0\leq \alpha\leq 1{\Big\}},\\ \mbox{if~~}0\neq x\in \mathbb{R}^m,\\
0, \mbox{otherwise}. \end{array}\right.\end{eqnarray*}

The operator $\mathrm{S}_{V_i}(x)$ was proposed in \cite{linren3}. Its role is to find the vector with the largest magnitude such that $\mathrm{S}_{V_i}(x)$ has the same direction as $x$, $\|\mathrm{S}_{V_i}(x)\|\leq\|x\|$ and $\alpha\mathrm{S}_{V_i}(x)\in V_i$ for all $0\leq \alpha\leq 1$.

 \begin{assumption}\label{ass2}{\rm \cite{linren3} Let $V_i\subseteq \mathbb{R}^m, i=1,\cdots,n$, be nonempty bounded closed
 sets such that $0\in V_i$, $\max_{x\in V_i}\|\mathrm{S}_{V_i}(x)\|=\bar{\rho}_i>0$ and $\inf_{x\notin V_i}\|\mathrm{S}_{V_i}(x)\|=\underline{\rho}_i>0$ for all $i$, where $\bar{\rho}_i$ and
 $\underline{\rho}_i$ are two positive
constants, and $\inf_{x\notin V_i}\|\mathrm{S}_{V_i}(x)\|$ denotes the infimum of $\|\mathrm{S}_{V_i}(x)\|$ when $x\notin V_i$.
}\end{assumption}

In Assumption \ref{ass2}, we do not require each $V_i$ to be convex. What we require on $V_i$ is that $V_i$ is bounded and the distance from any point outside $V_i$ to the origin is lower bounded by a positive constant.

\section{Problem Formulation and Algorithm}
Our objective is to design a distributed algorithm to make all agents cooperatively find the optimal state of
the optimization problem
\begin{eqnarray}\label{gels1}\begin{array}{lll}\mathrm{minimize}~~\sum_{i=1}^nf_i(s)\\
\mathrm{subject~to}~~s\in \mathbb{R}^m,\end{array}\end{eqnarray}
where $f_i(s): \mathbb{R}^m\rightarrow\mathbb{R}$ denotes the differentiable convex local objective
function that is known to only agent $i$.

\begin{assumption}\label{ass2u}{\rm \cite{linren4} Each set $X_i\triangleq{\Big\{}x{\Big |}\nabla f_i(x)=0{\Big\}}$, $i\in\mathcal{I}$, is nonempty and bounded.
%
}\end{assumption}

\begin{lemma}\label{lemma35}{\rm\cite{boyd} Let $f_0(\chi): \mathbb{R}^m\rightarrow\mathbb{R}$ be a differentiable convex function. $f_0(\chi)$ is minimized if and only if $\nabla f_0(\chi)=0$.}\end{lemma}
From Lemma \ref{lemma35}, each $X_i$ is the optimal set of the local objective function $f_i(x)$. Let $X$ be the optimal set of $\sum_{i=1}^nf_i(x)$. We have the following lemma.

\begin{lemma}\label{lemmallin}{\rm\cite{linren4} Under Assumptions \ref{ass2} and \ref{ass2u}, all $X_i$, $i\in\mathcal{I}$, and $X$ are nonempty closed bounded convex sets.}\end{lemma}

\begin{lemma}\label{lemma51}{\rm \cite{linren4} Let $f(s): \Xi\mapsto \mathbb{R}$ be a differentiable convex function and $Y$ be its minimum set in $\Xi$, where $\Xi\subseteq\mathbb{R}^m$ is a closed convex set. Suppose that $Y\subseteq \Xi$ is closed and bounded. For any $z=\lambda P_{Y}(y)+(1-\lambda)y$ with $\lambda\in(0,1)$, $0<\nabla f(z)^T\frac{y-P_{Y}(y)}{\|y-P_{Y}(y)\|}\leq \nabla f(y)^T\frac{y-P_{Y}(y)}{\|y-P_{Y}(y)\|}$ for any $y\in \Xi-Y$. }\end{lemma}

\begin{lemma}\label{lema18}{\rm \cite{linren4} Under Assumption \ref{ass2u}, $\lim_{\|y\|\rightarrow+\infty}f_i(y)=+\infty$ for all $i$ and accordingly $\lim_{\|y\|\rightarrow+\infty}\sum_{i=1}^nf_i(y)=+\infty$.}\end{lemma}

To solve the optimization problem (\ref{gels1}) in a distributed manner, we give the following algorithm for each agent by
\begin{eqnarray}\label{eq526}\begin{array}{lll}
{y}_i(k+1)={y}_i(k)+\mathrm{arctan} (e^{\|r_i(k)\|})T, y_i(0)>0,\\
\pi_i(k)=\sum\limits_{j\in \mathcal{N}_i(k)}a_{ij}(k)(r_j(k)-r_i(k))T,\\
q_i(k)=v_i(k)-p_i(k)v_i(k)T+\frac{p_i(k)}{2}\pi_i(k),\\
w_i(k)=r_i(k)+\frac{2}{p_i(k)}{v}_i(k)-v_i(k)T+\pi_i(k),\\
\theta_i(k)=\left\{\begin{array}{lll}
0, \mbox{~if~} \sqrt{y_i(k)}<\|\nabla
f_i(w_i(k))\|^2,\\
0, \mbox{~if~}q_i(k)-\frac{p_i(k)\nabla
f_i(w_i(k))}{2\sqrt{y_i(k)}}\\
\neq \mathrm{S}_{V_i}[q_i(k)-\frac{p_i(k)\nabla
f_i(w_i(k))}{2\sqrt{y_i(k)}}],\\
\frac{p_i(k)\nabla
f_i(w_i(k))}{2\sqrt{y_i(k)}}, \mbox{~otherwise,}\end{array}\right.\\
u_i(k)=\mathrm{S}_{V_i}[q_i(k)-\theta_i(k)],
\end{array}\end{eqnarray} for all $k\geq 0$,
where $v_i(0)\in {V_i}$ and $p_i(k)>0$ is the feedback damping gain.

In (\ref{eq526}),  $q_i(k)$ is a linear combination of the agent states, which is used to make all agents converge to a consensus point, and $\theta_i(k)$ contains a switching mechanism, where the first switching rule is used to guarantee all agents' position states to remain in a bounded region while the second switching rule is used to guarantee the balance of the optimal convergence rates of all agents. Specifically, in $\theta_i(k)$, the stepsize of the gradient, $\frac{1}{\sqrt{y_i(k)}}$, which is constructed based on only the position states, has two features: one is $\lim_{k\rightarrow+\infty}\frac{1}{\sqrt{y_i(k)}}=0$ and the other is $\lim_{k\rightarrow+\infty}\frac{\sqrt{y_j(k)}}{\sqrt{y_i(k)}}=1$ for all $i,j$, which will be shown later. The role of the stepsize $\frac{1}{\sqrt{y_i(k)}}$ is to minimize the effects of the gradient on the consensus convergence and keep the balance of the optimal convergence rates of all agents. $w_i(k)$ is a linear combination of the states and it is a key variable to determine the consensus behavior of the system dynamics, which will also be shown later.

In this note, our analysis is for the general $m$ case. When no confusion arises, the equations are given in the form of $m=1$ for simplicity of derivation expression.

\section{Stability analysis}
{Due to the nonconvexity of the operator $\mathrm{S}_{V_i}$ and the coexisting nonlinearities of the operator $\mathrm{S}_{V_i}$ and the gradient $\nabla f_i$, the system (\ref{eq11}) with (\ref{eq526}) is distinctly different from the distributed optimization problems with convex constraints and the consensus problem with nonconvex constraints studied in \cite{angelia, Elia, cotes, Cortes3, shi, liu, yuan, yuan1, Hong, Zhu, lup, Kvaternik, Zavlanos, srivast, plinwrenysong, linren4, linren5, linren3}. The approaches in \cite{angelia, Elia, cotes, Cortes3, shi, liu, yuan, yuan1, Hong, Zhu, lup, Kvaternik, Zavlanos, srivast, plinwrenysong, linren4, linren5, linren3} cannot be directly applied. }
To study the system (\ref{eq11}) with (\ref{eq526}), we first make some model transformations. {Let $\sigma_i(k)=\mathrm{S}_{V_i}[q_i(k)-\theta_i(k)]/[q_i(k)-\theta_i(k)].$
When $q_i(k)-\theta_i(k)=0$, $\mathrm{S}_{V_i}[q_i(k)-\theta_i(k)]=q_i(k)-\theta_i(k)$ and hence define
 $\sigma_i(k)=1$. Then $u_i(k)$ can be transformed into the form:   \begin{eqnarray*}\begin{array}{lll}u_i(k)=\mathrm{S}_{V_i}[q_i(k)-\theta_i(k)]=\sigma_i(k)[q_i(k)-\theta_i(k)].\end{array}\end{eqnarray*}}
As $\sigma_i(k)$ is a time-varying scaling factor, it is hard to perform analysis directly on the double integrator system with the control input in such a form. To proceed, we define two new variables $b_i(k)\in \mathbb{R}$ and $z_i(k)\in \mathbb{R}^m$ satisfying
 $1-b_i(k)T=\sigma_i(k)(1-p_i(k)T)$ and  $z_i(k)=r_i(k)+\frac{2}{p_i(k)}v_i(k)$. Note that \begin{eqnarray*}\begin{array}{lll}u_i(k)=\sigma_i(k)[v_i(k)-p_i(k)v_i(k)T+\frac{p_i(k)}{2}\pi_i(k)-\theta_i(k)]\\=v_i(k)-b_i(k)v_i(k)T+\frac{p_i(k)\sigma_i(k)}{2}\pi_i(k)-\sigma_i(k)\theta_i(k).\end{array}\end{eqnarray*}
Rewriting the system (\ref{eq11}) with (\ref{eq526}), we have
 \begin{eqnarray}\label{ef1s}\begin{array}{lll}r_i(k+1)=r_i(k)+\frac{p_i(k)}{2}[z_i(k)-r_i(k)]T\end{array}\end{eqnarray} and \begin{eqnarray}\label{ef2s}\begin{array}{lll}z_i(k+1)=z_i(k)+[\frac{p_i(k)T}{2}-1+\frac{p_i(k)}{p_i(k+1)}-\frac{p_i(k)b_i(k)T}{p_i(k+1)}]\\\times[z_i(k)-r_i(k)]+\frac{p_i(k)\sigma_i(k)}{p_i(k+1)}\pi_i(k)-\frac{2\sigma_i(k)}{p_i(k+1)}\theta_i(k).\end{array}\end{eqnarray} {It can easily be observed that the sums of the coefficients of $r_i(k)$ and $z_i(k)$ in (\ref{ef1s}) and (\ref{ef2s}) are both equal to $1$, which will be used for the system analysis.}
 It should be noted here that when $\theta_i(k)\neq0$, it follows from (\ref{eq526}) that $q_i(k)-\frac{p_i(k)\nabla
f_i(w_i(k))}{2\sqrt{y_i(k)}}
= \mathrm{S}_{V_i}[q_i(k)-\frac{p_i(k)\nabla
f_i(w_i(k))}{2\sqrt{y_i(k)}}]$ and hence $\sigma_i(k)=1$, i.e.,
  $b_i(k)=p_i(k)$. When $\theta_i(k)\neq0$, if $p_i(k+1)=b_i(k)$, we have that  \begin{eqnarray}\label{ef2x}\begin{array}{lll}w_i(k)=z_i(k)-\frac{p_i(k)}{2}[z_i(k)-r_i(k)]T+\pi_i(k).\end{array}\end{eqnarray}
  Rewriting (\ref{ef2s}) using $w_i(k)$ when $\theta_i(k)\neq0$ and $p_i(k+1)=b_i(k)$, we have that
  \begin{eqnarray}\label{ef2t}\begin{array}{lll}z_i(k+1)=w_i(k)-\frac{\nabla
f_i(w_i(k))}{\sqrt{y_i(k)}}.\end{array}\end{eqnarray}
Let $\phi(k)=[r_1(k)^T, z_1(k)^T,\cdots,r_n(k)^T, z_n(k)^T]^T$, $E_i(k)=\begin{bmatrix}1&0\\1&\frac{2}{p_i(k)}\end{bmatrix}$, $A_i(k)=\begin{bmatrix}1-\frac{p_i(k)T}{2}&\frac{p_i(k)T}{2}\\
{b_i(k)T}-\frac{p_i(k)T}{2}&1-{b_i(k)T}+\frac{p_i(k)T}{2}\end{bmatrix}$,
{$E(k)=\mathrm{diag}\{E_1(k),\cdots,E_n(k)\}$} and
$A(k)=\mathrm{diag}\{A_1(k),\cdots, A_n(k)\}$ where $E(k)$ and $A(k)$ are block diagonal matrices with their diagonal blocks equal to the matrices $E_i(k)$ and $A_i(k)$ respectively. Let
 {$\Lambda(k)=\mathrm{diag}\{\sigma_1(k),\cdots,\sigma_n(k)\}$, $W=\begin{bmatrix}0&0\\T&0\end{bmatrix}$} and $\nabla F(k)=[0^T,\frac{2\sigma_1(k)}{p_1(k+1)}\theta_1(k)^T,\cdots,0^T,\frac{2\sigma_n(k)}{p_n(k+1)}\theta_n(k)^T]^T$.

{It follows that
\begin{eqnarray}\label{eq23a}\hspace{-0.2cm}\begin{array}{lll}
&\phi(k+1)=\Psi(k)\phi(k)-\nabla F(k),\end{array}\end{eqnarray}}
\noindent where $\Psi(k)=E(k+1)E(k)^{-1}[A(k)-\Lambda(k)L(k)\otimes W]$ and $'\otimes'$ denotes the Kronecker product.

\begin{remark}{\rm
From the definition of $\sigma_i(k)$, each $\sigma_i(k)$ is time-varying and might not be uniform for all $i$ due to the nonconvex constraints, and hence the gradient weights might be nonuniform as well. From the view point of intuition, this
 might make the algorithm fail to solve the team optimization problem. Most of the existing works assumed the gradient weight to be uniform. Though articles \cite{linren4,linren5} have considered the case of nonuniform stepsizes, it is still unclear how to deal with the case in the system (\ref{eq23a}) when the nonuniform stepsizes and the time-varying scaling factors are taken into account simultaneously.}\end{remark}

In the following, we will first study the properties of the system matrices of (\ref{eq23a}) in Lemmas \ref{lemma1a}-\ref{lemma87}, and then the consensus and optimal convergence of (\ref{eq23a}) in Lemmas \ref{lemma1a1}-\ref{lemma9} and Theorem 1. Specifically, Lemma \ref{lemma1a} shows that the system matrices and the transition matrices are both stochastic. Lemma \ref{lemma1} shows that all agents' position states remain in a bounded region and each nonzero entry of these matrices is lower bounded by a positive constant by exploiting the swiching mechanism in the algorithm and the conditions given in the assumptions. Lemma \ref{lemma87} shows that there exist at least a column of the transition matrices each entry of which is lower bounded by a positive constant. Based on Lemma \ref{lemma87}, Lemma \ref{lemma1a1} shows that all agents reach a consensus as time evolves. Lemma \ref{lemma9} extends the continuous-time results of \cite{linren5} on gradient gains to the discrete-time system and  shows that the ratio of all stepsizes finally tends to $1$ as time evolves. Based on Lemmas \ref{lemma1a1} and \ref{lemma9}, Theorem 1 shows that the team optimization function is minimized as time evolves.
 \begin{assumption}\label{ass3}{\rm Suppose that $\frac{1}{T}>p_i(k+1)= b_i(k)>0$ for all $k\geq0$ and all $i$, and there exist a constant $d_{i}>0$ such that $p_i(k)> 2 d_{i}> 2 \lfloor {L}({k})\rfloor_{ii}$ for all $i$ and all $k\geq0$.}\end{assumption}

 Assumption \ref{ass3} actually gives a design rule for the algorithm, under which it will be shown that the transition matrices are all stochastic. Also, the constants $d_i$ always exist if $2\sum_{j\in \mathcal{N}_i(k)}a_{ij}(k)<p_i(0)$, which can be concluded from the proof of Lemma \ref{lemma1a}.

To implement the algorithm (\ref{eq526}) under Assumption  \ref{ass3}, we need to know the quantities,  $p_i(k)$, $w_i(k)$, $q_i(k)$, and $\theta_i(k)$. Since $p_i(k)=b_i(k-1)$, $p_i(k)$ can be obtained through computing $b_i(k-1)=[1-\sigma_i(k-1)(1-p_i(k-1)T)]/T$, where $\sigma_i(k-1)$ and $p_i(k-1)$ are both known at time instant $k$. In particular, for $k=0$, $p_i(0)$ can be adopted properly satisfying Assumption 3. Note that $w_i(k)$ and $q_i(k)$ are actually linear combinations of the variables $v_i(k)$, $r_i(k)$ and $r_j(k)$ for $j\in \mathcal{N}_i(k)$. Based on the obtained $p_i(k)$, $w_i(k)$ and $q_i(k)$ can be easily computed. As the variable $\theta_i(k)$ is dependent on the switching mechanism of the algorithm, we need judge the switching rules by computing $\sqrt{y_i(k)}$, $\|\nabla
f_i(w_i(k))\|^2$, $q_i(k)-\frac{p_i(k)\nabla
f_i(w_i(k))}{2\sqrt{y_i(k)}}$, and $\mathrm{S}_{V_i}[q_i(k)-\frac{p_i(k)\nabla
f_i(w_i(k))}{2\sqrt{y_i(k)}}]$ where $y_i(k)$ is known at time instant $k$ so as to obtain the variable $\theta_i(k)$. Though the algorithm computation looks a bit complex due to the existence of the switching mechanism, the algorithm does not require intermediate variables to be transmitted and it is a fully distributed algorithm.

Let {$\Gamma(k,s)=\prod_{i=s}^k\Psi(i)=\Psi(k)\cdots\Psi(s)$} be the transition matrix of the system (\ref{eq23a}).

 \begin{lemma}\label{lemma1a}{\rm Under Assumptions \ref{ass2}, \ref{ass2u} and \ref{ass3},
 $\Psi(k)$ and $\Gamma(k,s)$ are stochastic matrices for any $k\geq s\geq0$.}\end{lemma}
 \noindent{\textbf{Proof:}} {By simple calculations, 
\begin{eqnarray}\label{eq10013}E_i(k+1)E_i^{-1}(k)=\begin{bmatrix}1&0\\1-\frac{p_i(k)}{p_i(k+1)}&\frac{p_i(k)}{p_i(k+1)}\end{bmatrix}.\end{eqnarray}}
From the definition of $b_i(k)$, when $p_i(k)>0$, it is easy to see that $1>b_i(k)\geq p_i(k)$ and hence
\begin{eqnarray}\label{eq1001p} b_i(k+1)\geq p_i(k+1)\geq b_i(k)\geq p_i(k)>0\end{eqnarray} for all $k$ under Assumption \ref{ass3}. Thus $E(k+1)E(k)^{-1}$ is a stochastic matrix. Note that $L(k)\mathbf{1}=0$ from the definition of the graph Laplacian. It is easy to see that $[A(k)-\Lambda(k)L(k)\otimes W]\mathbf{1}=\mathbf{1}$. Note that $0\leq \sigma_i(k)\leq1$.
  Under Assumption \ref{ass3}, $b_i(k)T-\frac{p_i(k)T}{2}-\frac{p_i(k)\sigma_i(k)}{b_i(k)}d_iT\geq \frac{p_i(k)T}{2}-d_iT>0$ for all $i,k$.
 It can be easily checked that each entry of $A(k)-\Lambda(k)L(k)\otimes W$ is nonnegative and hence $A(k)-\Lambda(k)L(k)\otimes W$ is a stochastic matrix.
Therefore, $\Psi(k)$
 and $\Gamma(k,s)$ are both stochastic matrices for any $k\geq s\geq0$.\endproof

 {\begin{assumption}\label{ass13}{\rm  Suppose that there exist an infinite time sequence of $k_0,k_1,k_2,\cdots$ and a positive integer $\eta$ such that $k_0=0$, $0<k_{m+1}-k_m\leq \eta$ for all $m$ and the union of the graphs
$\mathcal{G}(k_m),\mathcal{G}(k_m+1),\cdots,\mathcal{G}(k_{m+1}-1)$
is strongly connected.}\end{assumption}}

Assumption \ref{ass13} ensures that the agents keep communication with each other persistently, which is a necessary condition for all agents to reach a consensus and minimize the team objective function.

\begin{lemma}\label{lemma1}{\rm Under Assumptions \ref{ass2}, \ref{ass2u}, \ref{ass3} and \ref{ass13},
\begin{itemize}\item[(1)] $\|r_i(k)\|<\rho$ and $\|z_i(k)\|<\rho$ for all $i,k$ and some constant $\rho>0$;
\item[(2)] each $\sigma_i(k)$ and each $\frac{p_i(k)}{p_i(k+1)}$ are both lower bounded by a positive constant for all $i,k$;
\item[(3)] each nonzero entry of  $A(k)-\Lambda(k)L(k)\otimes W$ is lower bounded by a positive constant.\end{itemize}
}\end{lemma}
\noindent{\textbf{Proof:}} Construct the Lyapunov function candidate $V(k)=\max_i\{\|r_i(k)-s\|,\|z_i(k)-s\|\}$ for some $s\in X$.
It is clear from Lemma \ref{lemma1a} that $\Psi(k)$ is stochastic. Hence, $\|r_i(k+1)-s\|\leq V(k)$ for all $i$.
In particular, when $\theta_i(k)=0$ for all $i$, $\|z_i(k+1)-s\|\leq V(k)$ also holds for all $i$ and thus $V(k+1)\leq V(k)$.

Suppose that $\theta_i(k)\neq0$. It is clear that $\sqrt{y_i(k)}\geq \|\nabla
f_i(w_i(k))\|^2$, and (\ref{ef2x})(\ref{ef2t}) hold under Assumption \ref{ass3}. It follows that 
\begin{eqnarray*}\begin{array}{lll}\|z_i(k+1)-s\|^2\\=\|w_i(k)-s\|^2+{\big\|}\frac{\nabla
f_i(w_i(k))}{\sqrt{y_i(k)}}{\big\|}^2-2(w_i(k)-s)^T\frac{\nabla
f_i(w_i(k))}{\sqrt{y_i(k)}}\\\leq\|w_i(k)-s\|^2+\frac{1}{\sqrt{y_i(k)}}-2(w_i-s)^T\frac{\nabla
f_i(w_i(k))}{\sqrt{y_i(k)}}.\end{array}\end{eqnarray*}
From the convexity of the function $f_i(\cdot)$, $-(w_i-s)^T\nabla
f_i(w_i(k))\leq f_i(s)-f_i(w_i(k))$. It follows that \begin{eqnarray}\label{eq0p0}\begin{array}{lll}\|z_i(k+1)-s\|^2\leq\|w_i(k)-s\|^2+\frac{1}{\sqrt{y_i(k)}}\\+2[f_i(s)-f_i(s_i)+f_i(s_i)-f_i(w_i(k))]/\sqrt{y_i(k)}\end{array}\end{eqnarray}
where $s_i\in X_i$. From the definition of $X_i$, $f_i(s)-f_i(s_i)\geq 0$ and $f_i(s_i)-f_i(w_i(k))\leq 0$. Note that $\pi/4\leq\mathrm{arctan}(e^{\|r_i(k)\|})\leq \pi/2$ for all $k$ and all $i$. There exists a constant $T_0>0$ such that $\frac{kT}{4}<y_i(k)<4kT$ for all $k\geq T_0$. \begin{eqnarray}\label{eq0p3}\begin{array}{lll}\|z_i(k+1)-s\|^2\leq\|w_i(k)-s\|^2+\frac{2}{\sqrt{kT}}\\+4[f_i(s)-f_i(s_i)]/\sqrt{kT}+[f_i(s_i)-f_i(w_i(k))]/\sqrt{kT}\end{array}\end{eqnarray}
for $k>T_0$.

Since $\lim_{\|y\|\rightarrow+\infty}f_i(y)=+\infty$ from Lemma \ref{lema18}, there exists one bounded convex closed region $\Omega_1=\{y\in \mathbb{R}^m \mid \|y-s\|\leq l_1\}$ for a constant $l_1>0$ such that  $f_i(w_i(k))-f_i(s)>4\max_j[f_j(s)-f_j(s_j)]+4n$ for any $w_i(k)\notin \Omega_1$.
Under Assumption \ref{ass3}, $0<p_i(k)T<1$ and hence $1-p_i(k)T/2\geq 1/2$.
If $w_i(k)\in \Omega_1$, then it follows from (\ref{ef2x}) and (\ref{eq0p3}) that $\|z_i(k)-s\|\leq 2l_1$ and 
$\|z_i(k+1)-s\|^2\leq l_1^2+2/\sqrt{T}+4[f_i(s)-f_i(s_i)]/\sqrt{T}\leq l_1^2+2/\sqrt{T}+4\max_i[f_i(s)-f_i(s_i)]/\sqrt{T}$.
 If $\|z_i(k)-s\|>2l_1$, then we have $w_i(k)\notin \Omega_1$ and hence $f_{i}(w_{i}(k))-f_{i}(s)\geq 4\max_i[f_j(s)-f_j(s_j)]+4n$. It follows from (\ref{eq0p3}) that $\|z_i(k+1)-s\|\leq \|w_i(k)-s\|\leq V(k)$.
Summarizing the above analysis, all agents remain in a bounded region. That is, $\|r_i(k)\|<\rho$ and $\|z_i(k)\|<\rho$ for all $i,k$ and some constant $\rho>0$. Hence, from the definitions of $z_i(k)$ and $\sigma_i(k)$, it can be obtained that $\|v_i(k)\|$ is bounded and hence from the definition of $\sigma_i(k)$, $\sigma_i(k)$ is lower bounded by a positive constant.

From the definition of $b_i(k)$, when $0<p_i(k)<\frac{1}{T}$, then $p_i(k)\leq b_i(k)$.
Under Assumption \ref{ass3}, $p_i(k)\leq b_i(k)\leq p_i(k+1)\leq b_i(k+1)$ and $p_i(k)T<1$. Hence, $1-p_i(k)T/2\geq 1/2$, $p_i(k)T/2\geq p_i(0)T/2$, $1-{b_i(k)T}+\frac{p_i(k)T}{2}\geq \frac{p_i(0)T}{2}$, and $\frac{p_i(k)}{p_i(k+1)}\geq \frac{p_i(0)}{p_i(k+1)}\geq p_i(0)T$. Consider the matrix
$A(k)-\Lambda(k)L(k)\otimes W$. ${b_i(k)T}-\frac{p_i(k)T}{2}-\sigma_i(k) d_iT\geq \frac{p_i(k)T}{2}-\sigma_i(k) d_iT\geq\frac{p_i(0)T}{2}-d_iT>0$.
Thus, each nonzero entry of  $A(k)-\Lambda(k)L(k)\otimes W$ is lower bounded by a positive constant.
\endproof

Theorem 3 in \cite{linren5} and Lemma 1 in \cite{linren3} both have considered the boundedness of the system states as well but both approaches there cannot be directly applied here because of the different adoption of the interaction mechanism or the lack of the consideration of the time-varying parameters.

\begin{lemma}\label{lemma87}{\upshape Under Assumptions \ref{ass2}, \ref{ass2u}, \ref{ass3} and \ref{ass13},
 there exists  a positive integer $h$ and a number  $0<\hat{\mu}<1$ such that
$\lfloor\Gamma(k_{m+4n}-1,k_{m})\rfloor_{ih}\geq \hat{\mu}$
  for all $k_m\geq0$ and $i$.
 }\end{lemma}
 \noindent{\textbf{Proof:}} The proof is very similar to the proof of Lemma 2 in \cite{linren3} and hence omitted. It should be noted that each agent might be the root node since the union of the graphs in $[k_m,k_{m+1}-1)$ is strongly connected under assumption 4. \endproof

\begin{lemma}\label{lemma1a1}{\rm Under Assumptions \ref{ass2}, \ref{ass2u}, \ref{ass3} and \ref{ass13}, $\lim_{m\rightarrow+\infty}[\phi_i(k)-\phi_j(k)]=0$ for all $i,j$, where $\phi_i(k),\phi_j(k)$ denote the $i,j$th entries of $\phi(k)$. }\end{lemma}
\noindent{\textbf{Proof:}} Since $\phi(k+1)=\Gamma(k,k)\phi(k)-\nabla F(k)$, then \begin{eqnarray}\label{efe3}\begin{array}{lll}\phi(k_{m+1})=\Gamma(k_{m+1}-1,k_m)\phi(k_m)\\-\sum_{j=k_m+1}^{k_{m+1}-1}\Gamma(k_{m+1}-1,j)\nabla F(j-1)-\nabla F(k_{m+1}-1).\end{array}\end{eqnarray} Note that $\frac{\pi}{4}\leq \mathrm{arctan} (e^{\|r_i(k)\|})\leq\frac{\pi}{2}$. There exists a constant $T_0>0$ such that $\frac{kT}{4}<y_i(k)<4kT$ for all $k>T_0$. Since all $\|r_i(k)\|$ and $\|z_i(k)\|$ are bounded from Lemma \ref{lemma1} and each $f_i(w_i(k))$ is differentiable and convex, each $\|w_i(k)\|$ is bounded and hence each $\|\nabla f_i(w_i(k))\|$ is bounded for all $i$. Thus, $\lim_{k\rightarrow+\infty}\|\nabla F(k)\|=0.$ There exists a constant $T_1>T_0$ for any $\epsilon>0$ such that $\|\nabla F(k)\|<\epsilon$ for all $k>T_1$. Since each $\Gamma(k,k)$ is a stochastic matrix from Lemma \ref{lemma1a} and $k_{m+1}-1-k_m\leq\eta$ under Assumption \ref{ass13}, it follows that \begin{eqnarray}\label{efe1}\begin{array}{lll}\|\sum_{j=k_m+1}^{k_{m+1}-1}\Gamma(k_{m+1}-1,j)\nabla F(j-1)\\-\nabla F(k_{m+1}-1)\|\leq (\eta+1)\epsilon\end{array}\end{eqnarray} for $k_m>T_1$. From Lemma \ref{lemma87}, there exists  a positive integer $i_0$ and a number  $0<\hat{\mu}< 1$ such that
$\lfloor\Gamma(k_{m+4n}-1,k_{m})\rfloor_{ii_0}\geq \hat{\mu}$
  for all $k_m\geq0$ and all $i$.
Note that \begin{eqnarray}\label{efe2}\begin{array}{lll}\|\phi_i(k_{m+4n})-\sum_{j\neq i_0}\lfloor\Gamma(k_{m+4n}-1,k_m)\rfloor_{ij}\phi_j(k_m)\\+\lfloor\Gamma(k_{m+4n}-1,k_m)\rfloor_{ii_0}\phi_{i_0}(k_m)\|\leq{4n(\eta+1)\epsilon}.\end{array}\end{eqnarray}It follows from (\ref{efe3}), (\ref{efe1}) and (\ref{efe2}) that \begin{eqnarray*}\hspace{-0.2cm}\begin{array}{lll}\max_i\phi_i(k_{m+4n})\leq (1-\hat{\mu})\max_i\phi_i(k_{m})+\hat{\mu}\phi_{i_0}(k_m)+4n(\eta+1)\epsilon,\\\min_i\phi_i(k_{m+4n})\geq (1-\hat{\mu})\min_i\phi_i(k_{m})+\hat{\mu}\phi_{i_0}(k_m)-4n(\eta+1)\epsilon,\end{array}\end{eqnarray*} and hence
\begin{eqnarray*}\begin{array}{lll}\max_i\phi_i(k_{m+4n})-\min_i\phi_i(k_{m+4n})\\\leq (1-\hat{\mu})[\max_i\phi_i(k_{m})-\min_i\phi_i(k_{m})]+8n(\eta+1)\epsilon.\end{array}\end{eqnarray*}
When $\max_i\phi_i(k_{m})-\min_i\phi_i(k_{m})\leq8n(\eta+1)\epsilon/\hat{\mu}$,
\begin{eqnarray*}\begin{array}{lll}\max_i\phi_i(k_{m+4n})-\min_i\phi_i(k_{m+4n})\\\leq [(1-\hat{\mu})/\hat{\mu}+1]8n(\eta+1)\epsilon\leq 8n(\eta+1)\epsilon/\hat{\mu}\leq 16n(\eta+1)\epsilon/\hat{\mu}.\end{array}\end{eqnarray*}
When $\max_i\phi_i(k_{m})-\min_i\phi_i(k_{m})>8n(\eta+1)\epsilon/\hat{\mu}$,
\begin{eqnarray*}\label{efe4}\begin{array}{lll}\max_i\phi_i(k_{m+4n})-\min_i\phi_i(k_{m+4n})\\-[\max_i\phi_i(k_{m})-\min_i\phi_i(k_{m})]<0.\end{array}\end{eqnarray*} This means that if $\max_i\phi_i(k_{m})-\min_i\phi_i(k_{m})\leq8n(\eta+1)\epsilon/\hat{\mu}$ holds, then
$\max_i\phi_i(k_{m+4ni})-\min_i\phi_i(k_{m+4ni})\leq8n(\eta+1)\epsilon/\hat{\mu}$ for all positive integers $i$.
Moreover, note that when $\max_i\phi_i(k_{m})-\min_i\phi_i(k_{m})>16n(\eta+1)\epsilon/\hat{\mu}$, \begin{eqnarray*}\label{efe4}\begin{array}{lll}\max_i\phi_i(k_{m+4n})-\min_i\phi_i(k_{m+4n})\\-[\max_i\phi_i(k_{m})-\min_i\phi_i(k_{m})]<-8n(\eta+1)\epsilon.\end{array}\end{eqnarray*} This means that all $\phi_i(k_{m+4n})$ will converge to the region in finite time where \begin{eqnarray*}\label{efe4}\begin{array}{lll}\max_i\phi_i(k_{m+4n})-\min_i\phi_i(k_{m+4n})\leq16n(\eta+1)\epsilon/\hat{\mu}.\end{array}\end{eqnarray*}

In view of the arbitrariness of the adoption of $k_m$ and $\epsilon$, letting $\epsilon\rightarrow0$, it follows that
$\lim_{m\rightarrow+\infty}[\phi_i(k_{m})-\phi_j(k_{m})]=0$ for all $i,j$. Note that $\|\phi(k+1)-\phi(k_m)\|\leq \|\Gamma(k,k_m)\phi(k_m)-\phi(k_m)\|+4n(\eta+1)\epsilon$ for all $k_m<k\leq k_{m+4n}-1$ and $\Gamma(k,k_m)$ is stochastic. Letting $k_m\rightarrow+\infty$ and $\epsilon\rightarrow0$,
it follows that $\lim_{k\rightarrow+\infty}\|\phi(k+1)-\phi(k_m)\|=\lim_{k\rightarrow+\infty}[\phi_i(k)-\phi_j(k)]=0$ for all $i,j$.\endproof

Lemma \ref{lemma1a1} is a key lemma to study the optimal convergence of the system (\ref{eq11}) with (\ref{eq526}). In contrast to Lemma 3 in \cite{linren3}, Lemma \ref{lemma1a1} need consider not only the interaction between agents but also the effects of the gradient term, which makes the analysis much more complicated.

\begin{lemma}\label{lemma9}{\rm For the system given by ${y}_i(k+1)={y}_i(k)+\mathrm{arctan} (e^{\|r_i(k)\|})T$ with $y_i(0)>0$, if $\lim_{k\rightarrow+\infty}[r_i(k)-r_j(k)]=0$ for all $i,j$, $\lim_{k\rightarrow+\infty}\frac{y_i(k)}{y_j(k)}=1$ for all $i,j$.}\end{lemma}

\noindent{\textbf{Proof:}} Note that $\pi/4\leq\mathrm{arctan}(e^{\|r_i(k)\|})\leq \pi/2$ for all $i$ and all $k$. There exists an integer $T_0>0$ such that $\frac{kT}{4}<y_i(k)<4kT$ for all $k>T_0$ and all $i$.
Let $\Delta_i(k)=y_i(0)-y_1(0)+\sum_{s=0}^{k}[\mathrm{arctan}(e^{\|r_i(s)\|})-\mathrm{arctan}(e^{\|r_1(s)\|})]T$. It is clear that $y_i(k)=y_1(k)(1+\Delta_i(k)/y_1(k))$ for all $i$. Since $\lim_{k\rightarrow+\infty}[r_i(k)-r_j(k)]=0$, from the continuity of the function $\mathrm{arctan}(e^{\|r_i(s)\|})$, there exists an integer $T_1>T_0$ for any ${\epsilon}_0>0$ such that
 $|\mathrm{arctan}(e^{\|r_i(k)\|})-\mathrm{arctan}(e^{\|r_1(k)\|})|<{\epsilon}_0$ for all $k>T_1$ and all $i$.
It is clear that $\lim_{k\rightarrow+\infty}\Delta_i(T_1)/y_1(k)=0$ and $|\sum_{s=T_1+1}^{k}[\mathrm{arctan}(e^{\|r_i(s)\|})-\mathrm{arctan}(e^{\|r_1(s)\|})]T/y_1(k)|$ $<4{\epsilon}_0$. Since ${\epsilon}_0$ can be arbitrarily chosen, we have
$\lim_{k\rightarrow+\infty}\Delta_i(k)/y_1(k)=\lim_{k\rightarrow+\infty}[\Delta_i(T_1)+\sum_{s=T_1+1}^{k}[\mathrm{arctan}(e^{\|x_i(s)\|})-\mathrm{arctan}(e^{\|x_1(s)\|})]]/y_1(k)=0$ and hence $\lim_{k\rightarrow+\infty}\frac{y_i(k)}{y_1(k)}=1$ for all $i$. Therefore, $\lim_{k\rightarrow+\infty}\frac{y_i(k)}{y_j(k)}=1$ for all $i,j$.\endproof

\begin{assumption}\label{ass2u1}{\rm Each graph $\mathcal{G}(k)$ is balanced, i.e., $\sum_{j=1}^na_{ij}(k)=\sum_{j=1}^na_{ji}(k)$.
%
}\end{assumption}

The role of Assumption \ref{ass2u1} is to balance the rate of the optimal convergence of all local objective functions.

\begin{theorem}\label{theorem1}{\rm Under Assumptions \ref{ass2}, \ref{ass2u}, \ref{ass3}, \ref{ass13} and \ref{ass2u1},
using (\ref{eq526}) for
(\ref{eq11}), all agents reach a
consensus and minimize the
 team objective function (\ref{gels1}) as $k\rightarrow+\infty$.
}
\end{theorem}
\noindent{\textbf{Proof:}} From Lemma \ref{lemma1a1}, $\lim_{k\rightarrow+\infty}[r_i(k)-z_j(k)]=0$ for all $i,j$. From the definition of $z_i(k)$, $\lim_{k\rightarrow+\infty}v_i(k)=0$ for all $i$. There exists a constant $T_0>0$ such that $\|v_i(k)\|<\underline{\rho}_i$ for all $i$ and all $k>T_0$. This means $\sigma_i(k)=1$ and hence $p_i(k)=b_i(k)$ and $p_i(k)=p_j(k+1)$ for all $i$ and $k>T_0$. As a result, the equations (\ref{ef1s}), (\ref{ef2x}) and (\ref{ef2t}) hold for $k>T_0$. Let $\phi^*(k)=\frac{1}{2n}\sum_{i=1}^{2n}\phi_i(k)$. Under Assumption \ref{ass2u1}, it can be checked that $\textbf{1}^T\Psi(k)=\textbf{1}^T$ for all $k$, where $\mathbf{1}$ denotes a column vector of all ones with a compatible dimension. It follows from (\ref{eq23a}) that \begin{eqnarray}\label{e1801}\begin{array}{lll}\phi^*(k+1)=\phi^*(k)-\frac{1}{2n}\mathbf{1}^T\nabla F(k).\end{array}\end{eqnarray}
Consider the Lyapunov function candidate $V_1(k)=\|\phi^*(k)-P_X(\phi^*(k))\|^2$ for $k>T_0$. From Lemma \ref{le1u}, $[\phi^*(k+1)-P_X(\phi^*(k+1))]^T[P_X(\phi^*(k))-P_X(\phi^*(k+1))]\leq0$. This implies that the angle between the vectors
$\phi^*(k+1)-P_X(\phi^*(k+1))$ and $P_X(\phi^*(k))-P_X(\phi^*(k+1))$ lies in $[\pi/2,\pi]$. From the triangle relationship, the angle between the vectors
$\phi^*(k+1)-P_X(\phi^*(k))$ and $P_X(\phi^*(k))-P_X(\phi^*(k+1))$ also lies in $[\pi/2,\pi]$. That is, $[\phi^*(k+1)-P_X(\phi^*(k))]^T[P_X(\phi^*(k))-P_X(\phi^*(k+1))]\leq0$. It follows from Lemma \ref{le1u} and (\ref{e1801}) that
\begin{eqnarray*}\hspace{-0.5cm}\begin{array}{lll}&&V_1(k+1)\\
&=&\|\phi^*(k+1)-P_X(\phi^*(k))+P_X(\phi^*(k))-P_X(\phi^*(k+1))\|^2\\
&=&\|\phi^*(k+1)-P_X(\phi^*(k))\|^2\\&+&\|P_X(\phi^*(k))-P_X(\phi^*(k+1))\|^2\\
&+&2[\phi^*(k+1)-P_X(\phi^*(k))]^T[P_X(\phi^*(k))-P_X(\phi^*(k+1))]\\
&\leq&\|\phi^*(k+1)-P_X(\phi^*(k))\|^2+\|\phi^*(k)-\phi^*(k+1)\|^2\\
&\leq&\|\phi^*(k+1)-P_X(\phi^*(k))\|^2+\|\frac{1}{2n}\mathbf{1}^T\nabla F(k)\|^2.
\end{array}
\end{eqnarray*}
From Lemma \ref{lemma1}, $\|r_i(k)\|$ and $\|z_i(k)\|$ are bounded for all $i,k$. Hence all $\|w_i(k)\|$ and $\phi^*(k)$ are bounded for all $i,k$. Since $f_i(w_i(k))$ is differentiable and convex, $\|\nabla f_i(w_i(k))\|$ is also bounded for all $i,k$. Note that $\frac{\pi}{4}\leq\mathrm{arctan} (e^{\|r_i(k)\|})\leq\frac{\pi}{2}$. There exists a constant $T_1>0$ such that $\frac{kT}{4}<y_i(k)<4kT$ for all $k>T_1$.
It follows from (\ref{e1801}) that
\begin{eqnarray*}\begin{array}{lll}&&V_1(k+1)\\
&\leq&\|\phi^*(k)-\frac{1}{2n}\mathbf{1}^T\nabla F(k)-P_X(\phi^*(k))\|^2+\|\frac{1}{2n}\mathbf{1}^T\nabla F(k)\|^2\\
&\leq&\|\phi^*(k)-P_X(\phi^*(k))\|^2-\frac{1}{n}(\phi^*(k)-P_X(\phi^*(k)))^T\\&\times&\sum_{i=1}^n\frac{\nabla f_i(w_i(k))^T}{\sqrt{y_i(k)}}+\|\frac{1}{n}\mathbf{1}^T\nabla F(k)\|^2\\
&\leq&\|\phi^*(k)-P_X(\phi^*(k))\|^2-\frac{1}{n}\sum_{i=1}^n(\phi^*(k)-P_X(\phi^*(k)))^T\\&\times&\frac{\nabla f_i(w_i(k))^T}{\sqrt{y_i(k)}}+\frac{\gamma }{k}\end{array}
\end{eqnarray*}for all $k>T_1$ and some constant $\gamma>0$.
 Note that $\lim_{k\rightarrow+\infty}[\phi_i(k)-\phi_j(k)]=0$ from Lemma \ref{lemma1a1}, $\frac{kT}{4}<y_i(k)<4kT$ for all $k>T_1$ and each function $f_i(\cdot)$ is differentiable. Then $\lim_{k\rightarrow+\infty}[w_i(k)-\phi^*(k)]=0$ for all $i$ and
there exists a constant $T_2>T_1$ for any $\epsilon_1>0$ such that $|w_i(k)-\phi^*(k)|<\epsilon_1$, $|\phi_i(k)-\phi^*(k)|<\epsilon_1$, and
$|f_i(w_i(k))-f_i(\phi^*(k))|<\epsilon_1$ for all $k>T_2$. From the definition of the projection operator, we have $\|\phi^*(k)-P_{X}(\phi^*(k))\|\leq \|\phi^*(k)-P_{X}(\phi^*(k-1))\|$.
Thus, \begin{eqnarray}\label{e1802}\begin{array}{lll}
&&V_1(k+1)-V_1(k)\\&\leq&
-\frac{1}{n}\sum_{i=1}^n(\phi^*(k)-w_i(k)+w_i(k)-P_X(\phi^*(k)))^T\\&\times&\frac{\nabla f_i(w_i(k))^T}{\sqrt{y_i(k)}}+\frac{\gamma }{k}\\
&\leq&-\frac{1}{n}\sum_{i=1}^n(w_i(k)-P_X(\phi^*(k)))^T\frac{\nabla f_i(w_i(k))^T}{\sqrt{y_i(k)}}\\
&+&\frac{\gamma }{k}+\frac{\tilde{\gamma}\epsilon_1}{\sqrt{k}}\\
&\leq&
-\frac{1}{n}\sum_{i=1}^n\frac{1}{\sqrt{y_i(k)}}[f_i(w_i(k))-f_i(P_{X}(\phi^*(k)))]\\&+&\frac{\gamma }{k}+\frac{\tilde{\gamma}\epsilon_1}{\sqrt{k}}\\
&\leq&
-\frac{1}{n}\sum_{i=1}^n\frac{1}{\sqrt{y_i(k)}}[f_i(w_i(k))-f_i(\phi^*(k))\\&+&f_i(\phi^*(k))-f_i(P_{X}(\phi^*(k)))]+
\frac{\gamma }{k}+\frac{\tilde{\gamma}\epsilon_1}{\sqrt{k}}\\
&\leq&-\frac{1}{n}\sum_{i=1}^n\frac{1}{\sqrt{y_i(k)}}[f_i(\phi^*(k))-f_i(P_{X}(\phi^*(k)))]\\&+&
\frac{\gamma }{k}+\frac{2\tilde{\gamma}\epsilon_1}{\sqrt{k}}
\end{array}\end{eqnarray} for $k>T_2$ and some constant $\tilde{\gamma}>0$,
where the third inequality has used the convexity of $f_i(w_i(k))$.

From Lemma \ref{lemma9}, $\lim_{k\rightarrow+\infty}\frac{y_i(k)}{y_j(k)}=1$ for all $i,j$. Together with the analysis above (\ref{e1802}), there must exist a constant $T_3>T_2$ such that
${\Big|}1-\sqrt{\frac{y_i(k)}{y_j(k)}}{\Big|}<\epsilon_1$ and ${\Big|}(\frac{\sqrt{y_1(k)}}{\sqrt{y_i(k)}}-1)[f_i(\phi^*(k))-f_i(P_{X}(\phi^*(k)))]{\Big|}<\epsilon_1$ for all $k>T_3$. Let $T_3$ be sufficiently large such that $\frac{1}{\sqrt{k}}<\epsilon_1$ for all $k>T_3$.
 It follows that  \begin{eqnarray*}\begin{array}{lll}
&&V_1(k+1)-V_1(k)\\&\leq&-\frac{1}{n}\sum_{i=1}^n\frac{1}{\sqrt{y_1(k)}}[f_i(\phi^*(k))-f_i(P_{X}(\phi^*(k)))]\\
&-&\frac{1}{n}\sum_{i=1}^n\frac{1}{\sqrt{y_1(k)}}(\frac{\sqrt{y_1(k)}}{\sqrt{y_i(k)}}-1)[f_i(\phi^*(k))\\&-&f_i(P_{X}(\phi^*(k)))]+
\frac{\gamma\epsilon_1}{\sqrt{k}}+\frac{2\tilde{\gamma}\epsilon_1}{\sqrt{k}}\\
&\leq&-\frac{1}{n\sqrt{y_1(k)}}[\sum_{i=1}^n(f_i(\phi^*(k))-f_i(P_{X}(\phi^*(k))))\\&-&n\epsilon_1]+\frac{\gamma\epsilon_1}{\sqrt{k}}+\frac{2\tilde{\gamma}\epsilon_1}{\sqrt{k}}\\
&\leq&-\frac{1}{2n\sqrt{kT}}[\sum_{i=1}^n(f_i(\phi^*(k))-f_i(P_{X}(\phi^*(k))))-\hat{\gamma}\epsilon_1]
\end{array}\end{eqnarray*}
 for all $k>T_3$ and some constant $\hat{\gamma}>0$.

{Let $\Omega_1=\{x\mid \|x-P_{X}(x)\|\leq l_1\}$ and $\Omega_2=\{x\mid \|x-P_{X}(x)\|\leq l_2\}$ for two constants $l_1>0$ and $l_2>0$ be two sets such that  $l_2>l_1+2\epsilon_1$ and $\min_{x\in \bar{\partial}\Omega_1}\sum_{i=1}^n[f_i(x)-f_i(s)]=\hat{\gamma}\epsilon_1+2\epsilon_1$ for any $s\in X$ where $\bar{\partial}\Omega_1$ denotes the boundary of $\Omega_1$. It is clear from Lemma \ref{lemma51} that $\sum_{i=1}^n(f_i(s_1)-f_i(s))>\hat{\gamma}\epsilon_1+2\epsilon_1$ for any $s_1\notin \Omega_1$.}
Note that each $\|\nabla f_i(w_i(k))\|$ is bounded. Let $T_3$ be further sufficiently large such that $\|\nabla F(k)\|<\epsilon_1$ for all $k>T_3$. From (\ref{e1801}), $\|\phi^*(k+1)-\phi^*(k)\|\leq \|\frac{1}{2n}\mathbf{1}^T\nabla F(k)\|<\epsilon_1$ for all $k>T_3$.
From the definition of $\Omega_2$, it is easy to see that for $k>T_3$, $\phi^*(k+1)\in \Omega_2$ when $\phi^*(k)\in \Omega_1$. When $\phi^*(k)\notin \Omega_1$ and $\phi^*(k)\in \Omega_2$, $V_1(k+1)-V_1(k)<0$ and hence $\|\phi^*(k+1)-P_{X}(\phi^*(k+1))\|\leq \|\phi^*(k)-P_{X}(\phi^*(k))\|$ for $k>T_3$, implying that
$\phi^*(k+1)\in \Omega_2$. By induction, it follows that if $\phi^*(\tilde{k})\in \Omega_2$ for some $\tilde{k}>T_3$, then $\phi^*({k})\in \Omega_2$ for all $k>\tilde{k}$.
When $\phi^*(k)\notin \Omega_2$ for $k>T_3$, $V_1(k+1)-V_1(k)<-\frac{1}{n\sqrt{kT}}\epsilon_1$.
Since $\sum_{k=h}^{+\infty}\frac{1}{\sqrt{kT}}=+\infty$ for some positive constant $h>0$, there exists a constant $T_4>T_3$ such that $\phi^*(k)\in \Omega_2$ for all $k>T_4$.
Since $\epsilon_1$ can be arbitrarily chosen, letting $\epsilon_1\rightarrow0$, it follows {from Lemma \ref{lemma51} and the definitions of $\Omega_1$ and $\Omega_2$} that $\lim_{k\rightarrow+\infty}\sum_{i=1}^n(f_i(\phi^*(k))-f_i(s))=0$. That is, the
 team objective function (\ref{gels1}) is minimized as $k\rightarrow+\infty$.\endproof

\section{Extension to the case with nonuniform convex position constraints}
In some applications, besides the velocity constraints, each agent's position state might also be constrained to a certain area.
In this section, we extend the results in previous sections to the case where there further exist nonuniform position constraints. {Here it is assumed that each agent's position state remains in a closed convex set, denoted by $H_i$, which is known to only agent $i$.} Under this circumstance, each agent has the following dynamics:
\begin{eqnarray}\label{eq112}\begin{array}{lll}
{r}_i(k+1)&=&P_{H_i}[r_i(k)+v_i(k)T]\\
{v}_i(k+1)&=&u_i(k),
\end{array}\end{eqnarray}
where $v_i$ is subject to the nonconvex constraint set $V_i$ as in previous sections. The Problem (\ref{gels1}) now becomes
\begin{eqnarray}\label{g00}\begin{array}{lll}\mathrm{minimize}~~\sum_{i=1}^nf_i(s)\\
\mathrm{subject~to}~~s\in H=\bigcap_{i=1}^nH_i.\end{array}\end{eqnarray}

To solve this propblem in a distributed manner, we propose the algorithm given by
\begin{eqnarray}\label{eq5216}\begin{array}{lll}
{y}_i(k+1)={y}_i(k)+\mathrm{arctan} (e^{\|r_i(k)\|})T, y_i(0)>0,\\

\pi_i(k)=\sum\limits_{j\in \mathcal{N}_i(k)}a_{ij}(k)(r_j(k)-r_i(k))T\\
q_i(k)=v_i(k)-p_iv_i(k)T+\frac{p_i}{4}\pi_i(k),\\
w_i(k)=r_i(k)+\frac{2}{p_i}{v}_i(k)-v_i(k)T+\frac{1}{2}\pi_i(k),\\
\theta_i(k)=\left\{\begin{array}{lll}
0, \mbox{~if~} S_{V_i}[q_i(k)-\frac{p_i\nabla
f_i(w_i(k))}{2\sqrt{y_i(k)}}]\\\neq q_i(k)-\frac{p_i\nabla
f_i(w_i(k))}{2\sqrt{y_i(k)}},\\
\frac{p_i\nabla
f_i(w_i(k))}{2\sqrt{y_i(k)}}, \mbox{~otherwise,}\end{array}\right.\\
u_i(k)=\mathrm{S}_{V_i}[q_i(k)-\theta_i(k)],
\end{array}\end{eqnarray} for all $k\geq 0$,
where $r_i(0)=P_{H_i}[r_i(0)]$, $v_i(0)=S_{V_i}(v_i(0))$ and $p_i>0$ is the feedback damping gain.  Here, the parameters $p_i$ are assumed to be constant for easy readability of our results.

 Compared with (\ref{eq526}), the switching mechanism in (\ref{eq5216}) is simplified. This is because all agents' states are bounded under Assumption \ref{ass2} and the the following assumption and there is no need to introduce switching rules to ensure the boundedness of the agents' states.

\begin{assumption}\label{as4}{\rm Each $H_i$ is closed and bounded for all $i$ and there exists a
scalar $\delta>0$ and a vector $\bar{x}\in H$ such that $\{\xi|\|\xi-\bar{x}\|\leq\delta\}\subset
H$.}\end{assumption}

Assumption \ref{as4} ensures that $H$ contains at least an interior point.

To present our main theorem under this situation, we need to modify Assumption \ref{ass3} as follows.

\begin{assumption}\label{ass31}{\rm Suppose that $\frac{1}{T}>p_i>0$ for all $i$, and there exist a constant $d_{i}>0$ such that $p_i> 2 d_{i}> 2 \lfloor {L}({k})\rfloor_{ii}$ for all $i$.}\end{assumption}

\begin{theorem}\label{theorem2}{\rm Under Assumptions \ref{ass2}, \ref{ass13}, \ref{ass2u1}, \ref{as4} and \ref{ass31},
using (\ref{eq5216}) for
(\ref{eq112}), all agents reach a
consensus and minimize the
 team objective function (\ref{gels1}) as $k\rightarrow+\infty$.
}
\end{theorem}
\noindent{\textbf{Sketch of Proof:}} Let $\sigma_i(k), z_i(k)$ and $b_i(k)$ be defined as previously. Specially, $z_i(k)=r_i(k)+\frac{2}{p_i}v_i(k)$. Under Assumption \ref{as4}, all $r_i(k)$ are bounded. From the definition of $z_i(k)$, under Assumption \ref{ass31}, it can be proved that all $z_i(k)$ are bounded and all $\sigma_i(k)$ are lower bounded by a positive constant. For some $i$, let $0\leq c_i(k)\leq 1$ be a scaling factor such that  $\varphi_{h}(k)\triangleq r_i(k)+c_i(k)\frac{p_i}{2}[z_i(k)-r_i(k)]T\in H_i$ and $r_i(k)+\frac{cp_i}{2}[z_i(k)-r_i(k)]T \notin H_i$ for all $1\geq c>c_i(k)$. Suppose that $c_i(k)<1$.
Let  $\varphi_{g}(k)=r_i(k)+\frac{p_i}{2}[z_i(k)-r_i(k)]T$, $\varphi_p(k)=P_{H_i}[\varphi_{g}(k)]$ and $\varphi_i(k)=\varphi_p(k)-\varphi_{h}(k)$.
Then by simple calculations similar to (\ref{ef2s}), we have \begin{eqnarray}\label{e890}\begin{array}{lll}z_i(k+1)=r_i(k+1)+\frac{2}{p_i}v_i(k+1)\\
=z_i(k)+[\frac{c_i(k)p_iT}{2}-b_i(k)T][z_i(k)-r_i(k)]+\varphi_i(k)\\
+\frac{\sigma_i(k)}{2}\pi_i(k)-\frac{2\sigma_i(k)}{p_i}\theta_i(k)\\
\triangleq c_z(k)z_i(k)+\sum_{j\in \{i\cup \mathcal{N}_i(k)\}}c_{jr}(k)r_j(k)\\+\varphi_i(k)-\frac{2\sigma_i(k)}{p_i}\theta_i(k)\\
\triangleq c_s(k)\varphi_{pz}(k)-\frac{2\sigma_i(k)}{p_i}\theta_i(k)\\+
(1-c_s(k))c_z(k)z_i(k)+\sum_{j\in \mathcal{N}_i(k)}c_{jr}(k)r_j(k).
\end{array}\end{eqnarray}
where $\varphi_{pz}(k)=\varphi_c(k)+\frac{1}{c_s}\varphi_i(k)$, $\varphi_c(k)=c_z(k)z_i(k)+\sum_{j\in \{i\cup \mathcal{N}_i(k)\}}c_{jr}(k)r_i(k)$,
and $c_s(k)=\frac{c_{ir}(k)}{\sum_{j\in \{i\cup \mathcal{N}_i(k)\}}c_{jr}(k)}$.
It is clear that $c_z(k)+\sum_{j\in \{i\cup \mathcal{N}_i(k)\}}c_{jr}(k)=c_s(k)+
(1-c_s(k))c_z(k)+\sum_{j\in \mathcal{N}_i(k)}c_{jr}(k)=1$ and $\sum_{j\in \{i\cup \mathcal{N}_i(k)\}}c_{jr}(k)=b_i(k)T-c_i(k)p_iT/2\geq b_i(k)T-p_iT/2$.
Under Assumption \ref{ass31}, $b_i(k)\geq p_i$ and $b_i(k)T<1$ as previously discussed. Hence $c_{ir}(k)\geq b_i(k)T-c_i(k)p_iT/2-d_iT/2\geq d_iT/2$ and $0<c_z(k)=1-b_i(k)T+c_i(k)p_iT/2\leq 1-b_i(k)T+p_iT/2\leq 1-p_iT/2$. Note that $c_{jr}(k)=a_{ij}(k)/2\geq0$. The coefficients of $z_i(k)$, $r_i(k)$ and $r_j(k)$ in (\ref{e890}) are all nonnegative. From the definition of the projection operator and Lemma 2 in \cite{yuan1}, $\|z_i(k+1)-P_{H}(z_i(k+1))\|\leq \|z_i(k+1)-P_{H}(z_i(k+1)+\frac{2\sigma_i(k)}{p_i}\theta_i(k))\|\leq \|z_i(k+1)+\frac{2\sigma_i(k)}{p_i}\theta_i(k)-P_{H}(z_i(k+1)+\frac{2\sigma_i(k)}{p_i}\theta_i(k))\|+\|\frac{2\sigma_i(k)}{p_i}\theta_i(k)\|\leq c_s(k)\|\varphi_{pz}(k)-P_{H}(\varphi_{pz}(k))\|+\|\frac{2\sigma_i(k)}{p_i}\theta_i(k)\|+
(1-c_s(k))c_z(k)\|z_i(k)-P_{H}(z_i(k))\|+\sum_{j\in \mathcal{N}_i(k)}c_{jr}(k)\|r_j(k)-P_{H}(r_j(k))\|$.

In the following, we study the term $\|\varphi_{pz}(k)-P_{H}(\varphi_{pz}(k))\|$ to analyze the convergence of $z_i(k)$.
Note that $\sum_{j\in \mathcal{N}_i(k)}c_{jr}(k)\leq d_iT/2$, $\sum_{j\in \{i\cup \mathcal{N}_i(k)\}}c_{jr}(k)\geq b_i(k)T-p_iT/2\geq p_iT/2$ and $p_i>2d_i$ under Assumption \ref{ass31}. Thus, $1/2<\frac{p_i-d_i}{p_i}\leq c_s(k)\leq1$. Let
$\varphi_z(k)=\varphi_c(k)+(1-c_i(k))(z_i(k)-r_i(k))p_iT/2/c_s(k)\triangleq z_i(k)(1-c_d(k))+r_i(k)c_d(k)$, where $c_{d}(k)=b_i(k)T-\frac{c_i(k)p_iT}{2}-(1-c_i(k))p_iT/2/c_s(k)$. It is clear that the triangles formed by the points $\varphi_{c}(k)$, $\varphi_{pz}(k)$ and $\varphi_{z}(k)$, and formed by $\varphi_{h}(k)$, $\varphi_{p}(k)$ and $\varphi_{g}(k)$ are similar. Moreover, note from (\ref{e890}) that $\varphi_c(k)=z_i(k)+[\frac{c_i(k)p_iT}{2}-b_i(k)T][z_i(k)-r_i(k)]$. Since $c_i(k)\leq 1$ and $b_i(k)\geq p_i$, then $b_i(k)T-\frac{c_i(k)p_iT}{2}\geq \frac{c_i(k)p_iT}{2}$ and hence $\varphi_c(k)\notin H_i$ based on the definition of $c_i(k)$ and the convexity of the convex set $H_i$. Let $\Gamma$ be a hyperplane such that $\varphi_z(k)-\varphi_{pz}(k)\perp\Gamma$ and $\varphi_p(k)\in \Gamma$. It is clear that $\varphi_g(k)-\varphi_{p}(k)\perp\Gamma$ and
all points of the convex set $H_i$ lie on one side of the hyperplane $\Gamma$. By considering the relationship between the aforementioned similar triangles, it can be obtained that $\varphi_z(k)$ and $\varphi_{pz}(k)$ lie on the other side of $\Gamma$. Recall that $\varphi_z(k)-\varphi_{pz}(k)\perp\Gamma$. The angle between the vectors $P_{H}(\varphi_z(k))-\varphi_{pz}(k)$ and $\varphi_z(k)-\varphi_{pz}(k)$ is no smaller than $\pi/2$. Thus, $\|\varphi_{pz}(k)-P_{H}(\varphi_z(k))\|\leq \|\varphi_z(k)-P_{H}(\varphi_z(k))\|$.
On the other hand,
 from Lemma 2 in \cite{yuan1}, $\|\varphi_{pz}(k)-P_{H}(\varphi_{pz}(k))\|\leq\|\varphi_{pz}(k)-P_{H}(\varphi_z(k))\|\leq \|\varphi_z(k)-P_{H}(\varphi_z(k))\|\leq (1-c_d(k))\|z_i(k)-P_{H}(z_i(k))\|+c_d(k)\|r_i(k)-P_{H}(r_i(k))\|$. It should be noted here that since
$1\leq 1/c_s\leq \frac{p_i}{p_i-d_i}<2$ and hence it can be proved that
$\rho_c< c_d(k)\leq b_i(k)T-p_iT/2$ for some constant $\rho_c>0$.

Summarizing the analysis above, it can be proved that $\|z_i(k+1)-P_{H}(z_i(k+1))\|\leq \tilde{c}_z(k)\|z_i(k)-P_{H}(z_i(k))\|+\|\frac{2\sigma_i(k)}{p_i}\theta_i(k)\|+\tilde{c}_{ir}(k)\|r_i(k)-P_{H}(r_i(k))\|+\sum_{j\in\mathcal{N}_i(k)}c_{jr}(k)\|r_j(k)-P_{H}(r_j(k))\|$ where $\tilde{c}_z(k)=c_s(k)(1-c_d(k))+(1-c_s(k))c_z(k)$, $\tilde{c}_{ir}(k)=c_s(k)c_d(k)$ and $c_{jr}(k)$ are lower bounded by a common positive constant, and $\tilde{c}_z(k)+\tilde{c}_{ir}(k)+\sum_{j\in\mathcal{N}_i(k)}c_{jr}(k)=1$ for all $j,k$. For the case of $c_i(k)=1$, similar statements can be obtained in a similar way. Following the lines of the proofs of Lemmas 9, 11, 12 and Proposition 1 in \cite{plinwrenysong}, it can be proved that $\lim_{k\rightarrow+\infty}[r_i(k)-z_j(k)]=0$  for all $i,j$. Further, by a similar analysis in the proof of Theorem \ref{theorem1}, it can be proved that the
 team objective function (\ref{g00}) is minimized as $k\rightarrow+\infty$.\endproof

\begin{remark}{\rm In our previous work \cite{linren4}, nonuniform stepsizes were considered but the analysis approach is hard to be applied directly here because this note is different in nature from \cite{linren4} in four aspects. First, the constraint sets considered in \cite{linren4} are uniform and convex while the constraint sets considered in this note are nonuniform and some might be nonconvex. Different convex constraint sets and nonconvex constraint sets might yield different nonlinearities. The coupling of different convex constraint sets and nonconvex constraint sets would yield more complicated nonlinearities. Second, the communication graphs in \cite{linren4} are assumed to be strongly connected at each time while the graphs in this note are assumed to be jointly strongly connected, which is more general and also much harder to analyze. Third, the agent dynamics is in the form of single integrators in \cite{linren4} while the agent dynamics is in the form of double integrators here in this note. Fourth, in \cite{linren4}, sign functions are used for the interactions between agents which can compensate for the inconsistent local gradients between neighbors and make the design and analysis relatively easier. In contrast, in this note, a kind of linear continuous consensus functions are used instead and no sign functions are used, which also makes the analysis of the system different from that in \cite{linren4}.}\end{remark}

 \section{Numerical Examples}
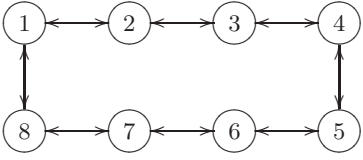
\begin{figure}
\begin{center}
\xymatrix{
    *++[o][F-]{1}\ar[r]\ar[d]
    & *++[o][F-]{2}\ar[r]\ar[l]
    & *++[o][F-]{3}\ar[r]\ar[l]
    & *++[o][F-]{4}\ar[d]\ar[l]
\\
    *++[o][F-]{8}\ar[u]\ar[r]
    & *++[o][F-]{7}\ar[l]\ar[r]
    & *++[o][F-]{6}\ar[l]\ar[r]
    & *++[o][F-]{5}\ar[l]\ar[u]
  }
\end{center}
\vspace{0.2cm}
\caption{One undirected graph.} \label{fig:1}
\end{figure}

\begin{figure}
\centering
\includegraphics[width=3.3in]{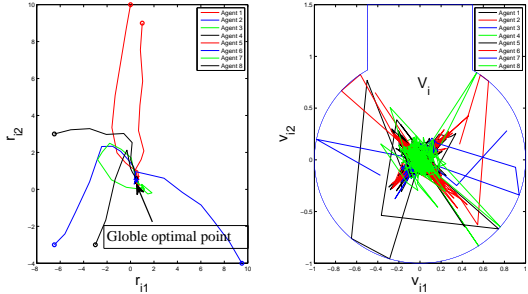} \\
\vspace{-0.4cm}
\caption{Trajectories of all agents with nonconvex velocity constraints and
nonuniform stepsizes}
\label{fig:3}
\end{figure}

\begin{figure}
\centering
\includegraphics[width=3.3in]{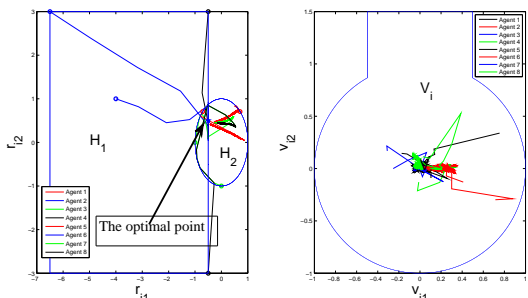} \\
\vspace{-0.4cm}
\caption{Trajectories of all agents with nonconvex velocity constraints, nonuniform position constraints and
nonuniform stepsizes}
\label{fig:4}
\end{figure}

Consider a multi-agent system with 8 agents in $\mathbb{R}^2$. {The communication graphs switch among the balanced subgraphs of the graph shown in Fig. \ref{fig:1}.}
 Each edge weight is 0.5. The
sample time is $T=0.2$ $s$ and the union of the communication graphs every $10$ $s$ is strongly connected.
The local objective functions for the agents are $f_1(r_1)=(r_{11}-1)^2+(r_{12}-1)^2$, $f_2(r_2)=r_{21}^2+r_{22}^2$, $f_3(r_3)=(r_{31}-1)^4+(r_{32}-1)^4$, $f_4(r_4)=r_{41}^4+r_{42}^4$, $f_5(r_5)=(r_{51}-1)^2+r_{52}^2$, $f_6(r_6)=r_{61}^2+(r_{62}-1)^2$, $f_7(r_7)=(r_{71}-1)^4+r_{72}^4$, and $f_8(r_8)=r_{81}^4+(r_{82}-1)^4$. The team performance function (\ref{gels1}) is minimized if and only if $r=[0.5,0.5]^T$.
 The velocity constraint set of each agent is $V_i=\{v\mid \|v\|\leq 1\}\cup \{v\mid -0.5\leq [1,0]^Tv\leq 0.5, 0\leq [0,1]^Tv\leq 1.5\}$ for all $i$. The position constraint set for agents $1,2,3,4$ are $H_1=\{r\mid \|r\|\leq 1\}$ while for agents $5,6,7,8$ is $H_2=\{r=[r_a,r_b]^T\mid -6.5\leq r_a\leq-0.5, -3\leq r_b\leq3 \}$.
 The team performance function (\ref{g00}) is minimized if and only if $r=[-0.5,0.5]^T$.
{According to Assumption \ref{ass3},  $p_i(0)$
is taken as $p_i(0)=1.5$ for algorithm (\ref{eq526}).} Fig. \ref{fig:3} shows the simulation results with nonconvex velocity constraints and
nonuniform stepsizes. It is clear that all agents reach a consensus and minimize the team performance function (\ref{gels1}) while their velocities remaining in $V_i$, which is consistent with Theorem \ref{theorem1}. According to Assumption \ref{ass31}, $p_i$ is taken as $p_i=1.5$ for algorithm (\ref{eq5216}).
 Fig. \ref{fig:4} show the simulation results with nonconvex velocity constraints, nonuniform position constraints and
nonuniform stepsizes. It is clear that all agents reach a consensus and minimize the team performance function (\ref{g00}) while their positions and velocities remaining in their corresponding constraint sets $H_i$ and $V_i$, which is consistent with Theorem \ref{theorem2} as well. 

\section{Conclusions}
In this note, a distributed optimization problem of multi-agent systems with nonconvex velocity constraints, nonuniform position constraints and
nonuniform stepsizes was studied.
 Two distributed constrained algorithms were proposed. The algorithm gains need not to be predesigned and can be given by each agent using its own and neighbours' information. The system considered was a nonlinear time-varying one and the analysis was performed based on a model transformation and the properties of stochastic matrices. It was shown that the optimization problem can be solved as long as the union of the communication graphs among each certain interval are jointly strongly connected and balanced.


\begin{thebibliography}{xx}


\bibitem{angelia} A. Nedi$\acute{\mathrm{c}}$, A. Ozdaglar, P. A. Parrilo, ``Constrained
consensus and optimization in multi-agent networks", IEEE
Transactions on Automatic Control, vol. 55, no. 4, pp.922-938, 2010.





\bibitem{Elia} J. Wang and N. Elia, ``A control perspective for centralized and distributed
convex optimization," in proceedings of IEEE Conference on Decision and Control, pp. 3800-3805, 2011.
\bibitem{cotes} B. Gharesifard and J. Cort$\mathrm{\acute{e}}$s, ``Distributed Continuous-Time Convex Optimization on Weight-Balanced Digraphs", IEEE
Transactions on Automatic Control, vol. 59, no. 3, pp. 781-786,
2014.

\bibitem{Cortes3}S. S. Kia, J. Cort$\acute{\mathrm{e}}$s, ``Distributed convex optimization via continuous-time coordination
algorithms with discrete-time communication", Automatica, vol. 55,  pp. 254-264, 2015.


\bibitem{shi} G. Shi, K. H. Johansson and Y. Hong, ``Reaching an Optimal Consensus: Dynamical Systems That Compute Intersections of Convex Sets", IEEE
Transactions on Automatic Control, vol. 58, no. 3, pp. 610-622,
2013.

\bibitem{liu} Z. Qiu, S. Liu and L. Xie, ``Distributed constrained optimal consensus of multi-agent systems", Automatica, vol. 68, pp. 209-216, 2016.








\bibitem{yuan} D. Yuan, Daniel W. C. Ho, Y. Hong, On convergence rate of distributed stochastic gradient algorithm for convex optimization with inequality constraints, SIAM Journal on Control and Optimization, vol. 54, no. 5, 2872-2892, 2016.

\bibitem{yuan1} P. Lin and  W. Ren, ``Constrained Consensus in Unbalanced Networks
With Communication Delays," IEEE Transations on Automatic Control, vol. 59, no. 3, pp. 775-781, 2014.
\bibitem{Hong} Y. Lou, Y. Hong and S. Wang, ``Distributed continuous-time approximate projection protocols for
shortest distance optimization problems", Automatica, vol. 69, no. 1, pp. 289-297, 2016.


\bibitem{Zhu} M. Zhu and S. Mart$\mathrm{\acute{{\mathrm{{{\imath}}}}}}$nez, ``An Approximate Dual Subgradient Algorithm for Multi-Agent Non-Convex Optimization", IEEE
Transactions on Automatic Control, vol. 58, no. 6, pp. 1534-1539.
2013.




\bibitem{lup}J. Lu, C. Y. Tang, ``Zero-gradient-sum algorithms for distributed convex
optimization: the continuous-time case", IEEE Transactions on Automatic Control,
vol. 57, no. 9, pp. 2348¨C2354, 2012.




\bibitem{Kvaternik}K. Kvaternik and L. Pavel, ``A Continuous-Time Decentralized Optimization Scheme With Positivity
Constraints", in proceedings of IEEE Conference on Decision and Control, pp. 6801-6807, 2012.



\bibitem{Zavlanos} N. Chatzipanagiotis and M. Zavlanos, ``A Distributed Algorithm for Convex Constrained
Optimization Under Noise", IEEE Transactions on Automatic Control,
vol. 61, no. 9, pp. 2496-2511, 2016.



\bibitem{srivast} K. Srivastava, and A. Nedi$\acute{\mathrm{c}}$, ``Distributed Asynchronous Constrained
Stochastic Optimization", IEEE Journal of Selected Topics in Signal
Processing, vol.5, no.4, pp.772-790,2011.

\bibitem{plinwrenysong} P. Lin,  W. Ren and Y. Song, ``Distributed Multi-agent Optimization Subject to Nonidentical Constraints and Communication Delays," Automatica, vol. 65, no.3, pp. 120-131, 2016.

    \bibitem{linren4} P. Lin,  W. Ren and J. A. Farrell, ``Distributed continuous-time optimization: nonuniform gradient gains, finite-time convergence, and convex constraint set," IEEE Transations on Automatic Control, vol. 62, no.5, pp. 2239-2253, 2017.

\bibitem{linren5} P. Lin, W. Ren, C. Yang and W. Gui. Distributed continuous-time and discrete-time optimization with nonuniform unbounded convex constraint sets and nonuniform stepsizes, IEEE Transactions on Automatic Control, accepted, DOI: 10.1109/TAC.2019.2910946.

\bibitem{linren3} P. Lin,  W. Ren and H. Gao, ``Distributed velocity-constrained consensus of discrete-time multi-agent
systems with nonconvex constraints, switching topologies, and delays," IEEE Transations on Automatic Control, vol. 62, no. 11, pp. 5788-5794, 2017.

\bibitem{s10} C. Godsil and G. Royle,
Algebraic Graph Theory. New York: Springer-Verlag, 2001.

\bibitem{boyd}S. Boyd and L. Vandenberghe. Convex Optimization, Cambridge University Press, 2004.
\bibitem{Facchinei}F. Facchinei and J. Pang. Finite-Dimensional Variational Inequalities and Complementarity Problems,  Springer-Verlag New York Inc., 2003.
\bibitem{s11} R. A. Horn and C. R. Johnson,
Matrix Analysis. Cambridge, U.K.: Cambridge Univ. Press, 1987.
\balance





















\end{thebibliography}
\end{document}